\documentclass[12pt,twoside]{article}
 \usepackage{amsmath}
 \usepackage{amssymb}
 \usepackage{array}
 \usepackage{graphicx}
 \usepackage{geometry}
\usepackage{graphics}
 \usepackage{enumerate}
 \usepackage[utf8]{inputenc}
\begin{document} \vskip3cm\noindent
 \setcounter{page}{1}
 \font\tinyfont=cmr8 \font\headd=cmr8
 \pagestyle{myheadings}
 \begin{center}
 \vskip.2cm\noindent {\bf Mellin Convolutions of Products and Ratios }\end{center}
 \vskip.5cm\noindent\begin{center} A.M. Mathai\\
 Department of Mathematics and Statistics,\\
McGill University, Montreal, Canada\\
 directorcms458@gmail.com\\
 and\\
 H.J. Haubold\\
 Office for Outer Space Affairs, United Nations\\
 Vienna International Centre, Vienna, Austria\\
 hans.haubold@gmail.com
 \end{center}

 \vskip.3cm\noindent{\bf Abstract}
 \vskip.3cm
 Usually, convolution refers to Laplace convolution in the literature. But Mellin convolutions can yield very ueful results. This aspect is illustrated in the coming sections. This paper deals with Mellin convolutions of products and ratios. Functions belonging to the pathway family of functions are considered. Several types of integral representations, their equivalent representations in terms of G and H-functions and their equivalent computable series representations are examined in this paper.



 \vskip.3cm\noindent{\bf 1.\hskip.3cm Introduction}
 \vskip.3cm
 Laplace convolutions involving products of two functions are widely known in the literature, see Luchko (2008) and Mainardi and Pagnini (2008). But Mellin convolutions of products and ratios are not widely used in the literature. Here we examine Mellin convolutions of products and ratios involving some functions belonging to the pathway family of functions. The pathway family of functions was introduced by Mathai (2005) for real-valued scalar functions with the arguments being rectangular matrices in the real domain. These results were extended to the complex domain in Mathai and Provost (2006). Here we consider the real scalar variable case first and then some generalizations are also considered. Let $x_1>0$ and $x_2>0$ be real positive scalar variables associated with the functions $f_1(x_1)$ and $f_2(x_2)$ respectively and with the joint function $f_1(x_1)f_2(x_2)$ the product. Then, the Mellin transforms of $f_1$ and $f_2$, with Mellin parameter $s$, and denoted by $M_{f_1}(s)$ and $M_{f_2}(s)$, are the following:

 $$M_{f_1}(s)=\int_0^{\infty}x_1^{s-1}f_1(x_1){\rm d}x_1\mbox{ and }M_{f_2}(s)=\int_0^{\infty}x_2^{s-1}f_2(x_2){\rm d}x_2\eqno(1.1)
 $$whenever the integrals are convergent. Let $u=x_1x_2$. Let $g_1(u)$ be the function associated with $u$. Then, the Mellin transform of $g_1(u)$ is the following:
 $$M_{g_1}(s)=\int_0^{\infty}u^{s-1}g_1(u){\rm d}u=\int_0^{\infty}\int_0^{\infty}(x_1x_2)^{s-1}f_1(x_1)f_2(x_2){\rm d}x_1\wedge{\rm d}x_2=M_{f_1}(s)M_{f_2}(s)\eqno(1.2)
 $$since the joint function of $x_1$ and $x_2$ is assumed to be the product $f_1(x_1)f_2(x_2)$. Consider the transformation
 $$u=x_1x_2,v=x_2\mbox{ or }v=x_1\Rightarrow {\rm d}x_1\wedge{\rm d}x_2=\frac{1}{v}{\rm d}u\wedge{\rm d}v
 $$and the joint function of $u$ and $v$ is $\frac{1}{v}f_1(\frac{u}{v})f_2(v)$ or $\frac{1}{v}f_1(v)f_2(\frac{u}{v})$. Then, $g_1(u)$ is the following:
 $$g_1(u)=\int_v\frac{1}{v}f_1(\frac{u}{v})f_2(v){\rm d}v=\int_v\frac{1}{v}f_1(v)f_2(\frac{u}{v}){\rm d}v.\eqno(1.3)
 $$But, from (1.2), $g_1(u)$ is available as the inverse Mellin transform of $M_{f_1}(s)M_{f_2}(s)$. That is,
 $$g_1(u)=\frac{1}{2\pi i}\int_{c-i\infty}^{c+i\infty}M_{f_1}(s)M_{f_2}(s)u^{-s}{\rm d}s,i=\sqrt{(-1)}\eqno(1.4)
 $$whenever the inverse Mellin transform exists. From (1.3) and (1.4) we have
 \begin{align*}
 g_1(u)&=\int_v\frac{1}{v}f_1(\frac{u}{v})f_2(v){\rm d}v=\int_v\frac{1}{v}f_1(v)f_2(\frac{u}{v}){\rm d}v\\
 &=\frac{1}{2\pi i}\int_{c-i\infty}^{c+i\infty}M_{f_1}(s)M_{f_2}(s)u^{-s}{\rm d}s\tag{1.5}\end{align*}
 This (1.5) is the Mellin convolution of a product property.
 \vskip.2cm
 When $x_1>0$ and $x_2>0$ are real scalar random variables with the densities $f_1(x_1)$ and $f_2(x_2)$ respectively and when $x_1$ and $x_2$ are statistically independently distributed, then
 \begin{align*}
 M_{f_1}(s)&=E[x_1^{s-1}], M_{f_2}(s)=E[x_2^{s-1}],\\
 E[(x_1x_2)^{s-1}]&=E[x_1^{s-1}]E[x_2^{s-1}]=M_{f_1}(s)M_{f_2}(s)\tag{1.6}\end{align*}
 due to statistical independence of $x_1$ and $x_2$, where $E[\cdot]$ denotes the expected value of $[\cdot]$. The density of $u=x_1x_2$, denoted by $g_1(u)$, is available as the first part of (1.5) if we use transformation of variables and as the second part of (1.5) if we use the result of arbitrary moments uniquely determining the corresponding density. Whenever $f_1(x_1)$ and $f_2(x_2)$ are statistical densities we can also give a statistical interpretation of the Mellin convolution of a product as the unique density of the product of independently distributed real scalar positive random variables $x_1$ and $x_2$. Thus, the results in (1.1) to (1.5) can be given statistical interpretations also whenever the functions $f_1$ and $f_2$ are statistical densities.

 \vskip.2cm
 For the real scalar positive variables case, the pathway family consists of the generalized type-1 beta family given by
 $$p_1(x)=\frac{\delta a^{\frac{\alpha}{\delta}}\Gamma(\frac{\alpha}{\delta}+\beta)}{\Gamma(\frac{\alpha}{\delta})\Gamma(\beta)}x^{\alpha-1}(1-ax^{\delta})^{\beta-1},0\le x\le 1,a>0,\alpha>0,\beta>0,1-ax^{\delta}>0\eqno(1.7)
 $$and zero elsewhere, with the standard form given by
 $$p_{11}(x)=\frac{\Gamma(\alpha+\beta)}{\Gamma(\alpha)\Gamma(\beta)}x^{\alpha-1}(1-x)^{\beta-1},0\le x\le 1,\alpha>0,\beta>0$$
 and zero elsewhere; the generalized type-2 beta model is given by
 $$p_2(x)=\frac{\delta a^{\frac{\alpha}{\delta}}\Gamma(\frac{\alpha}{\delta}+\beta)}{\Gamma(\frac{\alpha}{\delta})\Gamma(\beta)}x^{\alpha-1}(1+ax^{\delta})^{-(\frac{\alpha}{\delta}+\beta)},0\le x<\infty\eqno(1.8)
 $$for $a>0,\alpha>0,\beta>0,\delta>0$ and zero elsewhere, with the standard form given by
 $$p_{21}(x)=\frac{\Gamma(\alpha+\beta)}{\Gamma(\alpha)\Gamma(\beta)}x^{\alpha-1}(1+x)^{-(\alpha+\beta)},0\le x<\infty,\alpha>0,\beta>0$$
 and zero elsewhere; and the generalized gamma model is given by
 $$p_3(x)=\frac{\delta a^{\frac{\alpha}{\delta}}}{\Gamma(\frac{\alpha}{\delta})}x^{\alpha-1}{\rm e}^{-ax^{\delta}},0\le x<\infty,a>0,\delta>0,\alpha>0\eqno(1.9)
 $$and zero elsewhere with the standard form given by
 $$p_{31}(x)=\frac{a^{\alpha}}{\Gamma(\alpha)}x^{\alpha-1}{\rm e}^{-ax},0\le x<\infty,a>0,\alpha>0$$and zero elsewhere. In statistical densities the parameters are usually real and hence the conditions are stated for the real parameters but the functions are available when the parameters are in the complex domain also. In that case the same conditions hold on the real parts of the parameters. We will examine Mellin convolutions of products and ratios involving $f_1(x_1)$ and $f_2(x_2)$ belonging to the models (1.7) to (1.9) so that the results obtained will be readily applicable.
 \vskip.2cm
 The concept of Mellin convolution of a ratio is coming from the following considerations. As before, let $x_1>0$ and $x_2>0$ be real scalar positive variables with their associated functions $f_1(x_1)$ and $f_2(x_2)$ respectively and with the joint function $f_1(x_1)f_2(x_2)$ the product. Let $u=\frac{x_2}{x_1}$ the ratio and let $g_2(u)$ be the function corresponding to $u$. Then, proceeding as before we have
 $$E[u^{s-1}]=E[(\frac{x_2}{x_1})^{s-1}]=E[x_2^{s-1}]E[x_1^{-s+1}]=M_{f_2}(s)M_{f_1}(2-s).\eqno(1.10)
 $$Let $u=\frac{x_2}{x_1},v=x_2\Rightarrow {\rm d}x_1\wedge{\rm d}x_2=-\frac{v}{u^2}{\rm d}u\wedge{\rm d}v,x_2=v,x_1=\frac{v}{u}$. Also, $u=\frac{x_2}{x_1},v=x_1\Rightarrow {\rm d}x_1\wedge{\rm d}x_2=v{\rm d}u\wedge{\rm d}v,x_1=v,x_2=uv$. Then, the Mellin conovlution of a ratio property is the following:

 \begin{align*}
 g_2(u)&=\int_v\frac{v}{u^2}f_1(\frac{v}{u})f_2(v){\rm d}v=\int_vvf_1(v)f_2(uv){\rm d}v\\
 &=\frac{1}{2\pi i}\int_{c-i\infty}^{c+i\infty}M_{f_1}(2-s)M_{f_2}(s)u^{-s}{\rm d}s.\tag{1.11}\end{align*}Again, we will examine Mellin convolution of a ratio involving functions belonging to the pathway family (1.7) to (1.9) for the real scalar variable case first.
 \vskip.2cm
 This paper is organized as follows: Section 2 deals with Mellin convolution of a product involving functions in the standard form in (1.7) to (1.9). One general form is illustrated at the end. Section 3 gives Mellin convolutions of a ratio involving the standard forms in (1.7) to (1.9). Then, one illustration is given at the end for a general case. Section 4 examines some generalizations to the real matrix-variate case.

 \vskip.3cm\noindent{\bf 2.\hskip.3cm Mellin Convolution of a Product}
 \vskip.3cm
 Here we consider only Mellin convolutions of products of two real-valued scalar functions of real scalar positive variables. The results can be extended to Mellin convolution of a product involving three or more real-valued scalar functions.

 \vskip.3cm\noindent{\bf Problem 2.1.\hskip.3cm Type-2 beta versus gamma}
 \vskip.3cm To start with, we consider
 $$f_1(x_1)=\frac{\Gamma(\alpha+1+\beta)}{\Gamma(\alpha+1)\Gamma(\beta)}x_1^{\alpha}(1+x_1)^{-(\alpha+1+\beta)}, x_1\ge 0, f_2(x_2)=\frac{a^{\rho+1}}{\Gamma(\rho+1)}x_2^{\rho}{\rm e}^{-ax_2},x_2\ge 0,$$
 for $\alpha>-1,\rho>-1,\beta>0,a>0$ and let $f_1$ and $f_2$ be zero elsewhere. Let the joint function be $f_1(x_1)f_2(x_2)$ and $u=x_1x_2$ with the corresponding function $g_1(u)$. Then, for $c=\frac{\Gamma(\alpha+1+\beta)}{\Gamma(\alpha+1)\Gamma(\beta)}\frac{a^{\rho+1}}{\Gamma(\rho+1)}$,
 \begin{align*}
 g_1(u)&=\int_0^{\infty}\frac{1}{v}f_1(\frac{u}{v})f_2(v){\rm d}v
 =c\int_0^{\infty}\frac{1}{v}(\frac{u}{v})^{\alpha}(1+\frac{u}{v})^{-(\alpha+1+\beta)}v^{\rho}{\rm e}^{-av}{\rm d}v\\
 &=c~u^{\alpha}\int_0^{\infty}v^{\rho+\beta}(v+u)^{-(\alpha+1+\beta)}{\rm e}^{-av}{\rm d}v,0\le u<\infty\tag{2.1}\\
 g_1(u)&=\int_0^{\infty}\frac{1}{v}f_1(v)f_2(\frac{u}{v}){\rm d}v\\
 &=c~u^{\rho}\int_0^{\infty}v^{\alpha-\rho-1}(1+v)^{-(\alpha+1+\beta)}{\rm e}^{-\frac{au}{v}}{\rm d}v,0\le u<\infty.\tag{2.2}\end{align*}Also, the Mellin transforms are the following:

 \begin{align*}
 M_{f_1}(s)&=\frac{\Gamma(\alpha+1+\beta)}{\Gamma(\alpha+1)\Gamma(\beta)}\int_0^{\infty}x_1^{\alpha+s-1}(1+x_1)^{-(\alpha+1+\beta)}{\rm d}x_1\\
 &=\frac{\Gamma(\alpha+s)}{\Gamma(\alpha+1)}\frac{\Gamma(\beta-s+1)}{\Gamma(\beta)}, \Re(\alpha+s)>0,\Re(\beta-s+1)>0\tag{2.3}\\
 M_{f_2}(s)&=\frac{a^{\rho+1}}{\Gamma(\rho+1)}\int_0^{\infty}x_2^{\rho+s-1}{\rm e}^{-ax_2}{\rm d}x_2\\
 &=\frac{\Gamma(\rho+s)}{\Gamma(\rho+1)}\frac{a}{a^s}, \Re(\rho+s)>0\tag{2.4}\end{align*}where $\Re(\cdot)$ denotes the real part of $(\cdot)$. Then,
 \begin{align*}
 g_1(u)&=\frac{1}{\Gamma(\alpha+1)\Gamma(\beta)}\frac{a}{\Gamma(\rho+1)}\frac{1}{2\pi i}\int_{c-i\infty}^{c+i\infty}\Gamma(\alpha+s)\Gamma(\rho+s)\Gamma(\beta+1-s)(au)^{-s}{\rm d}s\tag{2.5}\\
 &=\frac{a}{\Gamma(\alpha+1)\Gamma(\beta)\Gamma(\rho+1)}G^{2,1}_{1,2}\left[au\vert_{\alpha.,\rho}^{-\beta}\right],0\le u<\infty\tag{2.6}\end{align*}where the $G(\cdot)$ is the G-function. For the theory and applications of the G-function, see for example, Mathai (1993). The inverse Mellin transform in (2.5) can be evaluated by using residue calculus. For $\alpha-\rho\ne\pm\nu,\nu=0,1,...$ the poles of the integrand are simple. Then, the function is available as the sum of the residues at the poles of $\Gamma(\alpha+s)$ and $\Gamma(\rho+s)$. These produce two confluent hypergeometric series and the sums of these residues is the explicit value of the G-function for $0\le u<\infty$. $g_1(u)$ above gives various representations of the density of $u$. Since the density is unique, we can compare the different representations of the unique density and obtain several results on the unique density or we can obtain several results on the G-function. Different types of  results on the G-function are given as Theorem 2.1 below. Since the procedure is the same, we will not give the derivations in detail for the remaining theorems in Section 2, only the results will be listed as theorems. The proofs are parallel to the steps given above.

 \vskip.3cm\noindent{\bf Theorem 2.1.}\hskip.3cm{\it For $a>0, u\ge 0,\alpha>-1,\rho>0,\beta>0$
 \begin{align*}
 G^{2,1}_{1,2}\left[au\vert_{\alpha,\rho}^{-\beta}\right]&=
 \Gamma(\alpha+1+\beta)a^{\rho}u^{\alpha}\int_{0}^{\infty}v^{\rho+\beta}(v+u)^{-(\alpha+1+\beta)}{\rm e}^{-av}{\rm d}v\\
 &=\Gamma(\alpha+1+\beta)a^{\rho}u^{\rho}\int_0^{\infty}v^{\alpha-\rho-1}(1+v)^{-(\alpha+1+\beta)}{\rm e}^{-a\frac{u}{v}}{\rm d}v\\
 &=(au)^{\alpha}\Gamma(\rho-\alpha)\Gamma(1+\alpha+\beta){_1F_1}(\beta+1+\alpha;1+\alpha-\rho:au)\\
 &+(au)^{\rho}\Gamma(1+\beta+\rho)\Gamma(\alpha-\rho){_1F_1}(\beta+1+\rho;1+\rho-\alpha;au),\end{align*}for $0\le au<\infty, \alpha-\rho\ne \pm\nu,n=0,1,...$.}

  \vskip.3cm\noindent Note that the integral can be evaluated at specific values of $a$ and $u$ such as $u=1$; $a=1$; $a=1,u=1$, giving rise to simpler integrals. The condition $\alpha-\rho\ne \pm\nu$ is needed only to write the G-function in terms of a confluent hypergeometric series. Each integral above multiplied by $\frac{1}{\Gamma(\alpha+1)\Gamma(\beta)}\frac{a}{\Gamma(\rho+1)}$ produces a statistical density $g_1(u)$ for $u$ also since our starting $f_1$ and $f_2$ are statistical densities. We have taken nonnegative integrable functions with total integral unity for convenience. The results hold for other types of functions provided the Mellin transforms exist and the product of the Mellin transforms is invertible. A direct generalization of Problem 2.1 can be achieved by replacing $1+x_1$ in $f_1$ by $1+a_1x_1^{\delta_1},a_1>0,\delta_1>0$ and replacing the exponent $ax_2$ in $f_2$ by $ax_2^{\delta_2},\delta_2>0$. Then, instead of the G-function, one will end up with a H-function and then one will obtain results on this H-function. For the theory and applications of the H-function, see for example, Mathai et al.(2010). Since $x_j\ge 0$ in Problem 2.1, we can obtain parallel results for $\delta_j<0$ also.

  \vskip.3cm\noindent{\bf Problem 2.2.\hskip.3cm Type-2 beta versus type-2 beta}
  \vskip.3cm
  Let
  $$f_j(x_j)=\frac{\Gamma(\alpha_j+1+\beta_j)}{\Gamma(\alpha_j+1)\Gamma(\beta_j)}x_j^{\alpha_j}(1+x_j)^{-(\alpha_j+1+\beta_j)},$$for $0\le x_j<\infty,\alpha_j>-1,\beta_j>0,j=1,2$
  and and let $f_j$ be zero elsewhere for $j=1,2$. Let the joint function of $x_1$ and $x_2$ be $f_1(x_1)f_2(x_2)$ and $u=x_1x_2$ with the corresponding function $g_1(u)$. Then, proceeding as in Problem 2.1 we have the following results:

 \vskip.3cm
 \noindent{\bf Theorem 2.2.}\hskip.3cm{\it For $\alpha_j>-1,\beta_j>0, j=1,2$ and for $c=\left\{\prod_{j=1}^2\frac{\Gamma(\alpha_j+1+\beta_j)}{\Gamma(\alpha_j+1)\Gamma(\beta_j)}\right\}$,
 \begin{align*}
 g_1(u)&=c~u^{\alpha_1}\int_0^{\infty}v^{\beta_1+\alpha_2}(u+v)^{-(\alpha_1+1+\beta_1)}(1+v)^{-(\alpha_2+1+\beta_2)}{\rm d}v,0\le u<\infty\\
 &=c~u^{\alpha_2}\int_0^{\infty}v^{\alpha_1+\beta_2}(1+v)^{-(\alpha_1+1+\beta_1)}(v+u)^{-(\alpha_2+1+\beta_2)}{\rm d}v, 0\le u<\infty\\
 &=\left\{\prod_{j=1}^2\frac{1}{\Gamma(\alpha_j+1)\Gamma(\beta_j)}\right\}G_{2,2}^{2,2}\left[u\vert_{\alpha_1,\alpha_2}^{-\beta_1,-\beta_2}\right],0\le u<\infty\end{align*}

 \begin{align*}
 G_{2,2}^{2,2}\left[u\vert_{\alpha_1,\alpha_2}^{-\beta_1,-\beta_2}\right]&=\Gamma(\alpha_2-\alpha_1)\Gamma(\beta_1+1+\alpha_1)\Gamma(\beta_2+1+\alpha_1)u^{\alpha_1}\\
 &\times {_2F_1}(\beta_1+1+\alpha_1,\beta_2+1+\alpha_1;1+\alpha_1-\alpha_2;u)\\
 &+\Gamma(\alpha_1-\alpha_2)\Gamma(\beta_1+1+\alpha_2)\Gamma(\beta_2+1+\alpha_2)u^{\alpha_2}\\
 &\times{_2F_1}(\beta_1+1+\alpha_2,\beta_2+1+\alpha_2;1+\alpha_2-\alpha_1;u)\\
 &\mbox{~~for }0\le u<1,\alpha_1-\alpha_2\ne\pm\nu,\nu=0,1,...\\
 &\Gamma(\alpha_1+1+\beta_1)\Gamma(\alpha_2+\beta_1+1)\Gamma(\beta_2-\beta_1)u^{-\beta_1-1}\\
 &\times{_2F_1}(\alpha_1+\beta_1+1,\alpha_2+\beta_1+1;1+\beta_1-\beta_2;\frac{1}{u})\\
 &+\Gamma(\alpha_1+\beta_2+1)\Gamma(\alpha_2+\beta_2+1)\Gamma(\beta_1-\beta_2)u^{-\beta_2-1}\\
 &\times{_2F_1}(\alpha_1+\beta_2+1,\alpha_2+\beta_2+1;1+\beta_2-\beta_1;\frac{1}{u})\\
 &\mbox{~~for }\beta_1-\beta_2\ne\pm\nu,\nu=0,1,...,1\le u<\infty\end{align*}}
 \vskip.3cm
 A direct generalization of Problem 2.2 is to replace $1+x_j$ in $f_j$ by $1+a_jx_j^{\delta_j},a_j>0,\delta_j>0,j=1,2$. Then, we will end up with a H-function, instead of the G-function above, and then parallel results will be available on this H-function. Since $x_j\ge 0$ in $f_j$ here, one can allow $\delta_j$ to be negative also. Parallel results will be available for $j=1,2$.

 \vskip.3cm\noindent{\bf Problem 2.3.\hskip.3cm Type-1 beta versus type-1 beta}
 \vskip.3cm
 Let
 $$f_j(x_j)=\frac{\Gamma(\alpha_j+1+\beta_j)}{\Gamma(\alpha_j+1)\Gamma(\beta_j)}x_j^{\alpha_j}(1-x_j)^{\beta_j-1},0\le x_j\le 1,\alpha_j>-1,\beta_j>0
 $$and zero elsewhere for $j=1,2$. Let the joint function be $f_1(x_1)f_2(x_2)$ and let $u=x_1x_2$ with the corresponding function denoted as $g_1(u)$. Then, proceeding as in Problem 2.1, one has the following results, for $\alpha_j>-1,\beta_j>0,j=1,2,c=\prod_{j=1}^2\frac{\Gamma(\alpha_j+1+\beta_j)}{\Gamma(\alpha_j+1)\Gamma(\beta_j)}$:

 \vskip.3cm\noindent{\bf Theorem 2.3.}\hskip.3cm{\it For $\alpha_j>-1,\beta_j>0,j=1,2$ and for $0\le u\le 1$
 \begin{align*}
 G^{2,0}_{2,2}[u\vert_{\alpha_1,\alpha_2}^{\alpha_1+\beta_1,\alpha_2+\beta_2}]&=\frac{1}{\Gamma(\beta_1)\Gamma(\beta_2)}u^{\alpha_1}\int_{v=u}^1v^{\alpha_2-\alpha_1-\beta_1}(v-u)^{\beta_1-1}(1-v)^{\beta_2-1}{\rm d}v,0\le u\le 1\\
 &=\frac{1}{\Gamma(\beta_1)\Gamma(\beta_2)}u^{\alpha_2}\int_{v=u}^1v^{\alpha_1-\beta_2-\alpha_2}(v-u)^{\beta_2-1}(1-v)^{\beta_1-1}{\rm d}v, 0\le u\le 1\\
 &=\frac{\Gamma(\alpha_2-\alpha_1)}{\Gamma(\beta_1)\Gamma(\alpha_2+\beta_2-\alpha_1)}u^{\alpha_1}\\
 &\times{_2F_1}(1-\beta_1,1+\alpha_1-\alpha_2-\beta_2;1+\alpha_1-\alpha_2;u)\\
 &+\frac{\Gamma(\alpha_1-\alpha_2)}{\Gamma(\beta_2)\Gamma(\alpha_1+\beta_1-\alpha_2)}u^{\alpha_2}\\
 &\times{_2F_1}(1-\beta_2,1+\alpha_2-\alpha_1-\beta_1;1+\alpha_2-\alpha_1;u),0\le u\le 1\\
 &\mbox{~~for }\alpha_1-\alpha_2\ne\pm\nu,\nu=0,1,...,0\le u< 1\end{align*}}
 \vskip.3cm
 This G-function multiplied by $\left\{\prod_{j=1}^2\frac{\Gamma(\alpha_j+1+\beta_j)}{\Gamma(\alpha_j+1)}\right\}$ gives a statistical density $g_1(u)$ for $u$ also. A direct generalization is available by replacing $1-x_j$ in $f_j$ by $1-a_jx_j^{\delta_j}$ for $a_j>0,\delta_j>0$, with the condition $1-a_jx_j^{\delta_j}>0,j=1,2$. Then, one will end up with a H-function, instead of the G-function above, and parallel results will be available for this H-function. It may be observed that a constant multiple of this H-function is a statistical density for $u$ also.

 \vskip.3cm\noindent{\bf Problem 2.4.\hskip.3cm Type-1 beta versus gamma}
 \vskip.3cm Let
 $$f_1(x_1)=\frac{\Gamma(\alpha+1+\beta)}{\Gamma(\alpha+1)\Gamma(\beta)}x_1^{\alpha}(1-x_1)^{\beta-1},0\le x_1\le 1, f_2(x_2)=\frac{a^{\gamma+1}}{\Gamma(\gamma+1)}x_2^{\gamma}{\rm e}^{-ax_2},x_2\ge 0
 $$for $\alpha>-1,\beta>0,\gamma>-1,a>0$ and let $f_1$ and $f_2$ be zero elsewhere and let the joint function be $f_1(x_1)f_2(x_2)$ and let $u=x_1x_2$ with the corresponding function denoted as $g_1(u)$. Then, proceeding as in Problem 2.1, we have the following results:

 \vskip.3cm\noindent{\bf Theorem 2.4.}\hskip.3cm{\it For $\alpha>-1,\gamma>-1,a>0,\beta>0$
 \begin{align*}
 G_{1,2}^{2,0}\left[au\vert_{\alpha,\gamma}^{\alpha+\beta}\right]&=\frac{a^{\gamma}}{\Gamma(\beta)}u^{\alpha}\int_{v=u}^{\infty}v^{\gamma-\alpha-\beta}(v-u)^{\beta-1}{\rm e}^{-av}{\rm d}v,0\le u<\infty\\
 &=\frac{a^{\gamma}}{\Gamma(\beta)}u^{\gamma}\int_0^1v^{\alpha-\gamma-1}(1-v)^{\beta-1}{\rm e}^{-a\frac{u}{v}}{\rm d}v, 0\le u<\infty\\
 &=\frac{\Gamma(\gamma-\alpha)}{\Gamma(\beta)}(au)^{\alpha}{_1F_1}(1-\beta;1+\alpha-\gamma:-au)\\
 &+\frac{\Gamma(\alpha-\gamma)}{\Gamma(\alpha+\beta-\gamma)}(au)^{\gamma}{_1F_1}(1+\gamma-\alpha-\beta;1+\gamma-\alpha;-au), 0\le u<\infty\\
 &\mbox{~~for }\alpha-\gamma\ne\pm\nu,\nu=0,1,...,0\le u<\infty\end{align*}}
 \vskip.3cm
 Note that this G-function multiplied by $\frac{\Gamma(\alpha+1+\beta)}{\Gamma(\alpha+1)}\frac{a}{\Gamma(\gamma+1)}$ produces a statistical density $g_1(u)$ for $u$ also. A direct generalization is available by replacing $1-x_1$ in $f_1$ by $1-a_1x_1^{\delta_1}$ with $1-a_1x_1^{\delta_1}>0, a_1>0,\delta_1>0$ and replacing the exponent $ax_2$ in $f_2$ by $ax_2^{\delta_2}$ with $\delta_2>0$. Then, one can obtain parallel results on a H-function. Here, $\delta_2$ can be negative also so that the integrals and the corresponding H-function will exist. If $\delta_1$ is allowed to be negative, then it will create problems with the support there and the problem will be complicated.

 \vskip.3cm\noindent{\bf Problem 2.5.\hskip.3cm Type-1 beta versus type-2 beta}
 \vskip.3cm Let
 $$f_1(x_1)=\frac{\Gamma(\alpha+1+\beta)}{\Gamma(\alpha+1)\Gamma(\beta)}x_1^{\alpha}(1-x_1)^{\beta-1}, f_2(x_2)=\frac{\Gamma(\gamma+1+\delta)}{\Gamma(\gamma+1)\Gamma(\delta)}x_2^{\gamma}(1+x_2)^{-(\gamma+1+\delta)},
 $$for $0\le x_1\le 1,x_2\ge 0,  \alpha>-1,\gamma>-1,\delta>0,\beta>0$ and let $f_1$ and $f_2$ be zero elsewhere. Let the joint function be $f_1(x_1)f_2(x_2)$ and let $u=x_1x_2$ with the corresponding function denoted as $g_1(u)$. Then, proceeding as in Problem 2.1, we will have the following results:

 \vskip.3cm\noindent{\bf Theorem 2.5.}\hskip.3cm{\it For $\alpha>-1,\gamma>-1,\beta>0,\delta>0$
 \begin{align*}
 G^{2,1}_{2,2}\left[u\vert_{\alpha,\gamma}^{-\delta,\alpha+\beta}\right]&=\frac{\Gamma(\gamma+1+\delta)}{\Gamma(\beta)}u^{\alpha}\int_{v=u}^{\infty}v^{\gamma-\alpha-\beta}(v-u)^{\beta-1}(1+v)^{-(\gamma+1+\delta)}{\rm d}v,0\le u<\infty\\
 &=\frac{\Gamma(\gamma+1+\delta)}{\Gamma(\beta)}u^{\gamma}\int_{v=0}^1v^{\alpha+\delta}(1-v)^{\beta-1}(v+u)^{-(\gamma+1+\delta)}{\rm d}v, 0\le u<\infty\\
 &=\begin{cases}\frac{\Gamma(\gamma-\alpha)\Gamma(1+\delta+\alpha)}{\Gamma(\beta)}u^{\alpha}{_2F_1}(1+\delta+\alpha,1-\beta;1+\alpha-\gamma;-u)\\
 +\frac{\Gamma(\alpha-\gamma)\Gamma(1+\delta+\gamma)}{\Gamma(\alpha+\beta-\gamma)}u^{\gamma}{_2F_1}(1+\delta+\gamma.1+\gamma-\alpha-\beta;1+\gamma-\alpha;-u)\\
 \mbox{ for }0\le u<1, \alpha-\gamma\ne\pm\nu,\nu=0,1,...\\
 \frac{\Gamma(\alpha+1+\delta)\Gamma(\gamma+1+\delta)}{\Gamma(\alpha+\beta+1+\delta)}u^{-1-\delta}{_2F_1}(\alpha+1+\delta,\gamma+1+\delta;\alpha+\beta+1+\delta;-\frac{1}{u})\\
 \mbox{ for }1\le u<\infty.\end{cases}\end{align*}}
 \vskip.3cm
 This G-function multiplied by $\frac{\Gamma(\alpha+1+\beta)}{\Gamma(\alpha+1)}\frac{1}{\Gamma(\gamma+1)\Gamma(\delta)}$ produces a statistical density $g_1(u)$ for $u$ also. A direct generalization is available by replacing $1-x_1$ in $f_1$ by $1-a_1x_1^{\delta_1}$ with $a_1>0,\delta_1>0,1-a_1x_1^{\delta_1}>0$ and replacing $1+x_2$ in $f_2$ by $1+a_2x_2^{\delta_2}$ with $a_2>0,\delta_2>0$. Then, one will obtain parallel results on a H-function. Also, note that in this Problem 2.5, if the type-2 beta density is replaced by a general density then the first integral representation gives Erd\'elyi-Kober fractional integral of the second kind of order $\beta$ and parameter $\alpha$, see also Mathai and Haubold (2018) and Fortin et al. (2016, 2020).

 \vskip.3cm\noindent{\bf Problem 2.6.\hskip.3cm Gamma versus gamma}
 \vskip.3cm Let
 $$f_j(x_j)=\frac{a_j^{\alpha_j+1}}{\Gamma(\alpha_j+1)}x_j^{\alpha_j}{\rm e}^{-a_jx_j},x_j\ge 0,\alpha_j>-1,a_j>0,j=1,2
 $$and let the joint function be $f_1(x_1)f_2(x_2)$ and $u=x_1x_2$ with the corresponding function denoted as $g_1(u)$. Then, proceeding as in Problem 2.1, we will have the following results:

 \vskip.3cm\noindent{\bf Theorem 2.6.}\hskip.3cm{\it For $a_j>0, \alpha_j>-1, j=1,2$
 \begin{align*}
 G^{2,0}_{0,2}\left[a_1a_2u\vert_{\alpha_1,\alpha_2}\right]&=a_1^{\alpha_1}a_2^{\alpha_2}u^{\alpha_1}\int_{v=0}^{\infty}v^{\alpha_2-\alpha_1-1}{\rm e}^{-a_1\frac{u}{v}-a_2v}{\rm d}v, 0\le u<\infty\\
 &=a_1^{\alpha_1}a_2^{\alpha_2}u^{\alpha_2}\int_0^{\infty}v^{\alpha_1-\alpha_2-1}{\rm e}^{-a_1v-a_2\frac{u}{v}}{\rm d}v, 0\le u<\infty\\
 &=(a_1a_2u)^{\alpha_1}\Gamma(\alpha_2-\alpha_1){_0F_1}(~~;1+\alpha_1-\alpha_2;a_1a_2u)\\
 &+(a_1a_1u)^{\alpha_2}\Gamma(\alpha_1-\alpha_2){_0F_1}(~~;1+\alpha_2-\alpha_1;a_1a_2u),0\le u<\infty\\
 &\mbox{~~for } \alpha_1-\alpha_2\ne\pm\nu,\nu=0,1,...\end{align*}}
 \vskip.3cm
 A direct generalization is available by replacing the exponent $a_jx_j$ by $a_jx_j^{\delta_j},\delta_j>0,j=1,2$. Then, we will obtain parallel results on a H-function. Here, one can allow $\delta_j<0,j=1,2$ also, which will then produce parallel results on a H-function. Since we have started with statistical densities for $f_1$ and $f_2$, we may observe that a constant multiple of this G-function in Theorem 2.6 is a statistical density $g_1(u)$ for $u$ also. When $\delta_1=1$ we have Kr\"atzel integral and associated with it is the Kr\"atzel transform, see Mathai and Haubold (2020). When $\delta_1=\delta_2=1$, that is the case considered in Theorem 2.6, the integrand in the integral representations, normalized, gives inverse Gaussian density. The unconditional density in the Bayesian analysis when the prior density is gamma and when the scale parameter has a prior gamma distribution then the unconditional density has the structure in the general as well as in the particular case where $\delta_1=\delta_2=1$. In this particular case, one has the Bessel integral given by the G-function representation in the simple poles cases also.

 \vskip.2cm For the sake of illustration, one direct generalization of a Mellin convolution of a product will be derived here in detail. This will be listed as Problem 2.7.

 \vskip.3cm\noindent{\bf Problem 2.7.\hskip.3cm Generalized type-2 beta versus generalized type-2 beta}
 \vskip.3cm
 Let
 $$f_j(x_j)=\frac{\delta_ja_j^{\frac{\alpha_j+1}{\delta_j}}\Gamma(\frac{\alpha_j+1}{\delta_j}+\beta_j)}{\Gamma(\frac{\alpha_j+1}{\delta_j})\Gamma(\beta_j)}x_j^{\alpha_j}(1+a_jx_j^{\delta_j})^{-(\frac{\alpha_j+1}{\delta_j}+\beta_j)},$$
 for $x_j\ge 0,a_j>0,\delta_j>0,\alpha_j>-1,j=1,2$ and $f_j$ is zero otherwise for $j=1,2.$ Let the joint function be $f_1(x_1)f_2(x_2)$ and let $u=x_1x_2$ and let the function corresponding to $u$ be denoted as $g_1(u)$. Then,
 \begin{align*}
 E[u^{s-1}]&=E[x_1^{s-1}]E[x_2^{s-1}]\\
 E[x_j^{s-1}]&=\frac{\delta_ja_j^{\frac{\alpha_j+1}{\delta_j}}\Gamma(\frac{\alpha_j+1}{\delta_j}+\beta_j)}{\Gamma(\frac{\alpha_j+1}{\delta_j})\Gamma(\beta_j)}\int_0^{\infty}x_j^{\alpha_j+s-1}(1+a_jx_j^{\delta_j})^{-(\frac{\alpha_j+1}{\delta_j}+\beta_j)}{\rm d}x_j\\
 &=\frac{a_j^{\frac{1}{\delta_j}}}{a_j^{\frac{s}{\delta_j}}}\frac{\Gamma(\frac{\alpha_j+s}{\delta_j})}{\Gamma(\frac{\alpha_j+1}{\delta_j})}\frac{\Gamma(\beta_j-\frac{s-1}{\delta_j})}{\Gamma(\beta_j)}.\tag{2.7}\end{align*}For convenience, let us consider the case $\delta_1=\delta_2=\delta$. Then, for $c_1=\prod_{j=1}^2\frac{a_j^{\frac{1}{\delta}}}{\Gamma(\frac{\alpha_j+1}{\delta})\Gamma(\beta_j)}$ and from (2.7) by taking the inverse Mellin transform we have $g_1(u)$ as the following:

 \begin{align*}
 g_1(u)&=c_1\frac{1}{2\pi i}\int_{c-i\infty}^{c+i\infty}[\prod_{j=1}^2\Gamma(\frac{\alpha_j+s}{\delta})\Gamma(\beta_j-\frac{s-1}{\delta})][(a_1a_2)^{\frac{1}{\delta}}u]^{-s}{\rm d}s\\
 &=c_1H_{2,2}^{2,2}\left[(a_1a_2)^{\frac{1}{\delta}}u\big\vert_{(\frac{\alpha_j}{\delta},\frac{1}{\delta}),j=1,2}^{(1-\beta_j-\frac{1}{\delta},\frac{1}{\delta}),j=1,2}\right],0\le u<\infty.\tag{2.8}\end{align*}The inverse Mellin transform in (2.8) will be evaluated by using residue calculus. Note that for $\alpha_1-\alpha_2\ne\pm\nu,\nu=0,1,...$ the poles of $\Gamma(\frac{\alpha_1+s}{\delta})\Gamma(\frac{\alpha_2+s}{\delta})$ are simple. The poles of $\Gamma(\frac{\alpha_j+s}{\delta})$ are at $\frac{\alpha_j+s}{\delta}=-\nu,\nu=0,1,...\Rightarrow s=-\alpha_j-\delta\nu,\nu=0,1,...$. The residue at $s=-\alpha_1-\delta\nu$ coming from $\Gamma(\frac{\alpha_1+s}{\delta})$ is $\lim_{s\to -\alpha_1-\delta\nu}(s+\alpha_1+\delta\nu)\Gamma(\frac{\alpha_1+s}{\delta})=\delta\frac{(-1)^{\nu}}{\nu!}$. The rest of the steps are parallel to those in the derivations in Problem 2.1. For $0\le a_1a_2u^{\delta}<1$ we will end up with two Gauss' hypergeometric series and the continuation part for $1\le a_1a_1u^{\delta}<\infty$ will give another two Gauss' hypergeometric series. If we consider the transformation of variables then $g_1(u)$ has the following forms where $c=\prod_{j=1}^2\frac{\delta a_j^{\frac{\alpha_j+1}{\delta}}\Gamma(\frac{\alpha_j+1}{\delta}+\beta_j)}{\Gamma(\frac{\alpha_j+1}{\delta})\Gamma(\beta_j)}$:
 \begin{align*}
 g_1(u)&=c\int_v\frac{1}{v}(\frac{u}{v})^{\alpha_1}[1+a_1(\frac{u}{v})^{\delta}]^{-(\frac{\alpha_1+1}{\delta}+\beta_1)}v^{\alpha_2}(1+a_2v)^{-(\frac{\alpha_2+1}{\delta}+\beta_2)}{\rm d}v\\
 &=c~u^{\alpha_1}\int_0^{\infty}v^{\alpha_2-\alpha_1-1}[1+\frac{a_1u^{\delta}}{v^{\delta}}]^{-(\frac{\alpha_1+1}{\delta}+\beta_1)}(1+a_2v)^{-(\frac{\alpha_2+1}{\delta}+\beta_2)}{\rm d}v.\tag{2.9}\end{align*}The second form is available by interchanging $(\alpha_1,\beta_1)$ and $(\alpha_2,\beta_2)$ in (2.9). Note that
 $$g_1(u)=\delta^2\left\{\prod_{j=1}^2a_j^{\frac{\alpha_j}{\delta}}\Gamma(\frac{\alpha_j+1}{\delta}+\beta_j)\right\}\times\mbox{ (the H-function)}$$where the H-function has the following explicit forms:

 \begin{align*}
 H&=H^{2,2}_{2,2}\left[(a_1a_2)^{\frac{1}{\delta}}u\big\vert_{(\frac{\alpha_j}{\delta},\frac{1}{\delta}),j=1,2}^{(1-\beta_j-\frac{1}{\delta},\frac{1}{\delta}),j=1,2}\right]\\
 &=\begin{cases}\delta\Gamma(\frac{\alpha_2-\alpha_1}{\delta})\Gamma(\beta_1+\frac{\alpha_1+1}{\delta})\Gamma(\beta_2+\frac{\alpha_1+1}{\delta})[(a_1a_2)^{\frac{1}{\delta}}u]^{\alpha_1}\\
 \times {_2F_1}(\beta_1+\frac{\alpha_1+1}{\delta},\beta_2+\frac{\alpha_1+1}{\delta};1+\frac{\alpha_1-\alpha_2}{\delta};a_1a_1u^{\delta})\\
 +\delta\Gamma(\frac{\alpha_1-\alpha_2}{\delta})\Gamma(\beta_1+\frac{\alpha_2+1}{\delta})\Gamma(\beta_2+\frac{\alpha_2+1}{\delta})[(a_1a_2)^{\frac{1}{\delta}}u]^{\alpha_2}\\
 \times {_2F_1}(\beta_1+\frac{\alpha_2+1}{\delta},\beta_2+\frac{\alpha_2+1}{\delta};1+\frac{\alpha_2-\alpha_1}{\delta};a_1a_2u^{\delta}),\\
 \mbox{ for }0\le a_1a_2u^{\delta}<1\mbox{ and for }\alpha_1-\alpha_2\ne\pm\nu,\nu=0,1,...\\
 \delta\Gamma(\beta_1+\frac{\alpha_1+1}{\delta})\Gamma(\beta_1+\frac{\alpha_2+1}{\delta})\Gamma(\beta_2-\beta_1)[(a_1a_2)^{\frac{1}{\delta}}u]^{-(\beta_1\delta+1)}\\
 \times {_2F_1}(\beta_1+\frac{\alpha_1+1}{\delta},\beta_1+\frac{\alpha_2+1}{\delta};1+\beta_1-\beta_2;\frac{1}{a_1a_2u^{\delta}})\\
 +\delta\Gamma(\beta_2+\frac{\alpha_1+1}{\delta})\Gamma(\beta_2+\frac{\alpha_2+1}{\delta})\Gamma(\beta_1-\beta_2)[(a_1a_2)^{\frac{1}{\delta}}u]^{-(\beta_2\delta+1)}\\
 \times {_2F_1}(\beta_2+\frac{\alpha_1+1}{\delta},\beta_2+\frac{\alpha_2+1}{\delta};1+\beta_2-\beta_1;\frac{1}{a_1a_2u^{\delta}})\\
 \mbox{ for }a_1a_2u^{\delta}\ge 1\mbox{ and for }\beta_1-\beta_2\ne\pm\nu,\nu=0,1,...\end{cases}\end{align*}This completes the calculations and the representations of $g_1(u)$ in two different integral representations and in terms of a H-function with explicit computable series form representations.

 \vskip.3cm\noindent{\bf 3.\hskip.3cm Mellin Convolution of a Ratio}
 \vskip.3cm
 Here we examine a ratio. Let $x_1>0$ and $x_2>0$ be real scalar positive variables associated with the functions $f_1(x_1)$ and $f_2(x_2)$ respectively. Let the joint function be $f_1(x_1)f_2(x_2)$. Let $u=\frac{x_2}{x_1}$ the ratio. Let the Mellin transforms of $f_1$ and $f_2$ be $M_{f_1}(s)$ and $M_{f_2}(s)$ respectively with the Mellin parameter $s$. Let the function associated with the ratio $u$ be denoted by $g_2(u)$. Then, the Mellin transform of $g_2(u)$ is available from the joint function as $M_{f_2}(s)M_{f_1}(2-s)$ which can easily be seen if $x_1$ and $x_2$ are statistically independently distributed random variables with the densities $f_1(x_1)$ and $f_2(x_2)$ respectively. Then,
 $$E[u^{s-1}]=E[(\frac{x_2}{x_1})^{s-1}]=E[x_2^{s-1}]E[x_1^{-s+1}]=M_{f_2}(s)M_{f_1}(2-s)\eqno(3.1)
 $$whenever the Mellin transforms exist. Then, from the inverse Mellin transform, $g_2(u)$ is available as the following:
 $$g_2(u)=\frac{1}{2\pi i}\int_{c-i\infty}^{c+i\infty}M_{f_2}(s)M_{f_1}(2-s)u^{-s}{\rm d}s,i=\sqrt{(-1)}\eqno(3.2)
 $$whenever the integral is convergent. We can also reach $g_2(u)$ through transformation of variables. Let $x_2=v$ and $u=\frac{x_2}{x_1}$ then $x_1=\frac{v}{u}$ and ${\rm d}x_1\wedge{\rm d}x_2=-\frac{v}{u^2}{\rm d}u\wedge{\rm d}v$. If one takes $x_1=v$ then $x_2=uv$ and ${\rm d}x_1\wedge{\rm d}x_2=v{\rm d}u\wedge{\rm d}v$. Then, we have two different integral representations for $g_2(u)$, namely
 \begin{align*}
 g_2(u)&=\int_v\frac{v}{u^2}f_1(\frac{v}{u})f_2(v){\rm d}v\tag{3.3}\\
 &=\int_vvf_1(v)f_2(uv){\rm d}v.\tag{3.4}\end{align*}
 Then, (3.2),(3.3) and (3.4) give three different representations for the same function $g_2(u)$. We will consider various functions belonging to the pathway family of functions for convenience.

 \vskip.3cm\noindent{\bf Problem 3.1.\hskip.3cm Gamma versus gamma}
 \vskip.3cm
 Let
 $$f_j(x_j)=\frac{a_j^{\alpha_j+1}}{\Gamma(\alpha_j+1)}x_j^{\alpha_j}{\rm e}^{-a_jx_j},x_j\ge 0,\alpha_j>-1,a_j>0,j=1,2
 $$and let the joint function be $f_1(x_1)f_2(x_2)$ and let $c=\prod_{j=1}^2\frac{a_j^{\alpha_j+1}}{\Gamma(\alpha_j+1)}$. Then, from (3.3)
 $$g_2(u)=c\int\frac{v}{u^2}(\frac{v}{u})^{\alpha_1}{\rm e}^{-a_1\frac{v}{u}}v^{\alpha_2}{\rm e}^{-a_2v}{\rm d}v=c\Gamma(\alpha_1+\alpha_2+2)u^{\alpha_2}(a_1+a_2u)^{-(\alpha_1+\alpha_2+2)}\eqno(3.5)
 $$and from (3.4)
 $$g_2(u)=c\int_vvv^{\alpha_1}{\rm e}^{-a_1v}(uv)^{\alpha_2}{\rm e}^{-a_2uv}{\rm d}v=c~u^{\alpha_2}\Gamma(\alpha_1+\alpha_2+2)(a_1+a_2u)^{-(\alpha_1+\alpha_2+2)},\eqno(3.6)
 $$for $0\le u<\infty,$ which  is the same as the one in (3.5). Also,
 $$M_{f_2}(s)=\frac{a_2^{\alpha_2+1}}{\Gamma(\alpha_2+1)}\int_0^{\infty}x_2^{\alpha_2+s-1}{\rm e}^{-a_2x_2}{\rm d}x_2=\frac{a_2}{a_2^s}\frac{\Gamma(\alpha_2+s)}{\Gamma(\alpha_2+1)}, \Re(\alpha_2+s)>0\eqno(3.7)
 $$and
 $$M_{f_1}(2-s)=\frac{a_1^{\alpha_1+1}}{\Gamma(\alpha_1+1)}\int_0^{\infty}x_1^{\alpha_1-s+1}{\rm e}^{-a_1x_1}{\rm d}x_1=\frac{a_1^s}{\Gamma(\alpha_1+1)a_1}\Gamma(2+\alpha_1-s),\Re(2+\alpha_1-s)>0.\eqno(3.8)
 $$Then,
 \begin{align*}
 g_2(u)&=\frac{a_2}{a_1}\frac{1}{\Gamma(\alpha_1+1)\Gamma(\alpha_2+1)}\frac{1}{2\pi i}\int_{c-i\infty}^{c+i\infty}\Gamma(\alpha_2+s)\Gamma(2+\alpha_1-s)(\frac{a_2u}{a_1})^{-s}{\rm d}s\\
 &=\frac{a_2}{a_1\Gamma(\alpha_1+1)\Gamma(\alpha_2+1)}G^{1,1}_{1,1}\left[\frac{a_2u}{a_1}\vert_{\alpha_2}^{-1-\alpha_1}\right],0\le u<\infty.\tag{3.9}\end{align*}This G-function can be evaluated as the sum of the residues at the poles of $\Gamma(\alpha_2+s)$ for $0\le \frac{a_2u}{a_1}<1$ and as the sum of the residues at the poles of $\Gamma(2+\alpha_1-s)$ for $1\le \frac{a_2u}{a_1}<\infty$. Thus, we have

 $$g_2(u)=\frac{\Gamma(\alpha_1+\alpha_2+2)}{\Gamma(\alpha_1+1)\Gamma(\alpha_2+1)}\frac{1}{u}\begin{cases}(\frac{a_2u}{a_1})^{\alpha_2+1}{_1F_0}(2+\alpha_1+\alpha_2;~~;-\frac{a_2u}{a_1}),0\le \frac{a_2u}{a_1}<1\\
 (\frac{a_1}{a_2u})^{\alpha_1+1}{_1F_0}(2+\alpha_1+\alpha_2;~~;-\frac{a_1}{a_2u}),1\le \frac{a_2u}{a_1}<\infty.\end{cases}\eqno(3.10)$$Then, from (3.5),(3.6) and (3.10) we have the following theorem:

 \vskip.3cm\noindent{\bf Theorem 3.1.}\hskip.3cm{\it For $a_j>0,\alpha_j>-1, c=\left\{\prod_{j=1}^2\frac{a_j^{\alpha_j+1}}{\Gamma(\alpha_1+1)}\right\}$,
 \begin{align*}
 g_2(u)&=c\Gamma(\alpha_1+\alpha_2+2)u^{\alpha_2}(a_1+a_2u)^{-(\alpha_1+\alpha_2+2)},0\le u<\infty\\
 &=\frac{\Gamma(\alpha_1+\alpha_2+2)}{\Gamma(\alpha_1+1)\Gamma(\alpha_2+1)}\frac{1}{u}\begin{cases}(\frac{a_2u}{a_1})^{\alpha_2+1}{_1F_0}( 2+\alpha_1+\alpha_2;~~;-\frac{a_2u}{a_1}),0\le \frac{a_2u}{a_1}<1\\
 (\frac{a_1}{a_2u})^{\alpha_1+1}{_1F_0}(2+\alpha_1+\alpha_2;~~;-\frac{a_1}{a_2u}),1\le \frac{a_2u}{a_1}<\infty,\end{cases}\end{align*}}where $g_2(u)$ is a statistical density for $u$ also.
 \vskip.3cm
 Generalization is available by replacing $a_jx_j$ in $f_j(x_j)$ by $a_jx_j^{\delta_j},\delta_j>0,j=1,2.$ In this case, one will end up with a H-function instead of the G-function above. Then, various integral and series representations of that H-function will be available. Here, we will be able to consider $\delta_j<0,j=1,2$ also giving rise to parallel results. Since the derivations of the results in Problem 3.1 are given in detail the remaining results will be stated without the detailed derivations. The derivations will be parallel to those provided in Problem 3.1.

 \vskip.3cm\noindent{\bf Problem 3.2.\hskip.3cm Gamma versus type-2 beta}
 \vskip.3cm
 Let
 $$f_1(x_1)=\frac{a_1^{\alpha_1+1}}{\Gamma(\alpha_1+1)}x_1^{\alpha_1}{\rm e}^{-a_1x_1}, f_2(x_2)=\frac{\Gamma(\alpha_2+1+\beta_2)}{\Gamma(\alpha_2+1)\Gamma(\beta_2)}x_2^{\alpha_2}(1+x_2)^{-(\alpha_2+1+\beta_2)}
 $$for $\alpha_j>-1,0\le x_j<\infty,j=1,2,a_1>0,\beta_2>0$ and let the joint function be $f_1(x_1)f_2(x_2)$. Let $u=\frac{x_2}{x_1}$ with the corresponding function denoted as $g_2(u)$. Then, following through the derivations parallel to those in Problem 3.1 we have the following results:

 \vskip.3cm\noindent{\bf Theorem 3.2.}\hskip.3cm{\it For $\alpha_j>-1, j=1,2,\beta_2>0, a_1>0$
 \begin{align*}
 &G^{1,2}_{2,1}\left[\frac{u}{a_1}\vert_{\alpha_2}^{-\beta_2,-1-\alpha_1}\right]\\
 &=a_1^{\alpha_1+2}\Gamma(\alpha_2+\beta_2+1)u^{-\alpha_1-2}\int_0^{\infty}v^{\alpha_1+\alpha_2+1}(1+v)^{-(1+\alpha_2+\beta_2)}{\rm e}^{-a_1\frac{v}{u}}{\rm d}v,0\le u<\infty\\
 &=a_1^{\alpha_1+2}\Gamma(\alpha_2+\beta_2+1)u^{\alpha_2}\int_0^{\infty}v^{\alpha_1+\alpha_2+1}(1+uv)^{-(1+\alpha_2+\beta_2)}{\rm e}^{-a_1v}{\rm d}v,0\le u<\infty\\
 &=\Gamma(\alpha_2+\beta_2+1)\Gamma(1+\alpha_1-\beta_2)(\frac{a_1}{u})^{1+\beta_2}{_1F_1}(1+\alpha_2+\beta_2;\beta_2-\alpha_1;\frac{a_1}{u})\\
 &+\Gamma(2+\alpha_1+\alpha_2)\Gamma(\beta_2-\alpha_1-1)(\frac{a_1}{u})^{2+\alpha_1}{_1F_1}(2+\alpha_1+\alpha_2;2+\alpha_1-\beta_2;\frac{a_1}{u}),0\le u<\infty\\
 &\mbox{~~for  }\beta_2-\alpha_1-1\ne\pm\nu,\nu=0,1,...\end{align*}}

 \vskip.3cm
 Generalizations are possible by replacing the exponent $a_1x_1$ in $f_1$ by $a_1x_1^{\delta_1},\delta_1>0$ and replacing $1+x_2$ in $f_2$ by $1+a_2x_2^{\delta_2},a_2>0,\delta_2>0$. Then, one will end up with a H-function and results will be available for this H-function. In this case, $\delta_j<0,j=1,2$ will also work and parallel results will be available.

 \vskip.3cm\noindent{\bf Problem 3.3.\hskip.3cm Type-2 beta versus gamma}
 \vskip.3cm
 Let
 $$f_1(x_1)=\frac{\Gamma(\gamma+1+\delta)}{\Gamma(\gamma+1)\Gamma(\delta)}x_1^{\gamma}(1+x_1)^{-(\gamma+1+\delta)},f_2(x_2)=\frac{a^{\alpha+1}}{\Gamma(\alpha+1)}x_2^{\alpha}{\rm e}^{-ax_2}
 $$for $\gamma>-1,\alpha>-1,a>0,\delta>0,0\le x_j<\infty,j=1,2$ and let the joint function be $f_1(x_1)f_2(x_2)$. Let $u=\frac{x_2}{x_1}$. Then, proceeding as in the derivations in Problem 3.1, we will have the following results:

 \vskip.3cm\noindent{\bf Theorem 3.3.}\hskip.3cm{\it For $\alpha>-1,\gamma>-1,\delta>0,a>0$
 \begin{align*}
 G^{2,1}_{1,2}\left[au\vert_{\alpha,\delta-1}^{-1-\gamma}\right]
 &=\Gamma(\gamma+1+\delta)a^{\alpha}u^{\delta-1}\int_0^{\infty}v^{\alpha+\gamma+1}(u+v)^{-(\gamma+1+\delta)}{\rm e}^{-av}{\rm d}v,0\le u<\infty\\
 &=\Gamma(\gamma+1+\delta)a^{\alpha}u^{\alpha}\int_0^{\infty}v^{\alpha+\gamma+1}(1+v)^{-(\gamma+1+\delta)}{\rm e}^{-auv}{\rm d}v,0\le u<\infty\\
 &=\Gamma(\delta-1-\alpha)\Gamma(2+\gamma+\alpha)(au)^{\alpha}{_1F_1}(2+\gamma+\alpha;2+\alpha-\delta;au)\\
 &+\Gamma(1+\gamma+\delta)\Gamma(1+\alpha-\delta)(au)^{\delta-1}{_1F_1}(1+\gamma+\delta;\delta-\alpha;au),0\le u<\infty\\
 &\mbox{~~for }\delta-1-\alpha\ne\pm\nu,\nu=0,1,...\end{align*}}

 \vskip.3cm
 Generalizations are possible by replacing $1+x_1$ in $f_1$ by $1+b_1x_1^{\delta_1},b_1>0,\delta_1>0$ and by replacing the exponent in $f_2$ by $ax_2^{\delta_2},\delta_2>0$. This will produce a H-function and several results on this H-function. In this case, $\delta_j<0,j=1,2$ will also work, resulting in parallel results.

 \vskip.3cm\noindent{\bf Problem 3.4.\hskip.3cm Gamma versus type-1 beta}
 \vskip.3cm Let
 $$f_1(x_1)=\frac{a^{\gamma+1}}{\Gamma(\gamma+1)}x_1^{\gamma}{\rm e}^{-ax_1}, f_2(x_2)=\frac{\Gamma(\alpha+1+\beta)}{\Gamma(\alpha+1)\Gamma(\beta)}x_2^{\alpha}(1-x)^{\beta-1}
 $$for $a>0,\beta>0,\gamma>-1,\alpha>-1, x_1\ge 0, 0\le x_2\le 1$ and $f_1$ and $f_2$ are zero otherwise. Let the joint function be $f_1(x_1)f_2(x_2)$ and let $u=\frac{x_2}{x_1}$. Then, proceeding as in the derivations for Problem 3.1 we have the following results:

 \vskip.3cm\noindent{\bf Theorem 3.4.}\hskip.3cm{\it For $a>0, \beta>0,\alpha>-1,\gamma>-1$
 \begin{align*}
 G^{1,1}_{2,1}\left[\frac{u}{a}\vert_{\alpha}^{-1-\gamma,\alpha+\beta}\right]&=\frac{a^{\gamma+2}}{\Gamma(\beta)}u^{-\gamma-2}\int_0^1v^{\alpha+\gamma+1}(1-v)^{\beta-1}{\rm e}^{-a\frac{v}{u}}{\rm d}v, 0<u<\infty\\
 &=\frac{a^{\gamma+2}}{\Gamma(\beta)}u^{\alpha}\int_0^{\frac{1}{u}}v^{\alpha+\gamma+1}(1-uv)^{\beta-1}{\rm e}^{-av}{\rm d}v, 0\le u<\infty\\
 &=\frac{\Gamma(2+\alpha+\gamma)}{\Gamma(2+\alpha+\gamma+\beta)}(\frac{u}{a})^{-\gamma-2}{_1F_1}(2+\alpha+\gamma;2+\alpha+\gamma+\beta;-\frac{a}{u}),\end{align*}}

 \vskip.3cm
 for $0\le u<\infty$. Generalization can be achieved by replacing the exponent $ax_1$ in $f_1$ by $ax_1^{\delta_1},\delta_1>0$ and replacing $1-x_2$ in $f_2$ by $1-a_2x_2^{\delta_2}, a_2>0,\delta_2>0, 1-a_2x_2^{\delta_2}>0$. Then, one will end up with a H-function and several results on this H-function.

 \vskip.3cm\noindent{\bf Problem 3.5.\hskip.3cm Type-1 beta versus gamma}
 \vskip.3cm Let
 $$f_1(x_1)=\frac{\Gamma(\alpha+1+\beta)}{\Gamma(\alpha+1)\Gamma(\beta)}x_1^{\alpha}(1-x_1)^{\beta-1}, f_2(x_2)=\frac{a^{\gamma+1}}{\Gamma(\gamma+1)}x_2^{\gamma}{\rm e}^{-ax_2}
 $$for $\alpha>-1,\gamma>-1,a>0,\beta>0,x_2\ge 0,0\le x_1\le 1$ and $f_1$ and $f_2$ are zero otherwise. Let the joint function be $f_1(x_1)f_2(x_2)$ and let $u=\frac{x_2}{x_1}$. Then, proceeding as in the derivations in Problem 3.1, we have the following results:

 \vskip.3cm\noindent{\bf Theorem 3.5.}\hskip.3cm{\it For $\alpha>-1,\gamma>-1,a>0,\beta>0$
 \begin{align*}
 G^{1,1}_{1,2}\left[au\vert_{\gamma,-1-\alpha-\beta}^{-1-\alpha}\right]&=\frac{a^{\gamma}}{\Gamma(\beta)}u^{-\beta-\alpha-1}\int_{v=0}^uv^{\alpha+\gamma+1}(u-v)^{\beta-1}{\rm e}^{-av}{\rm d}v,0<u<\infty\\
 &=\frac{a^{\gamma}}{\Gamma(\beta)}u^{\gamma}\int_0^1v^{\alpha+\gamma+1}(1-v)^{\beta-1}{\rm e}^{-auv}{\rm d}v,0<u<\infty\\
 &=a^{\gamma}\frac{\Gamma(2+\alpha+\gamma)}{\Gamma(2+\alpha+\gamma+\beta)}u^{\gamma}{_1F_1}(2+\alpha+\gamma;2+\alpha+\gamma+\beta;-au),0\le u<\infty.\end{align*}}

 \vskip.3cm A direct generalization is achieved by replacing $1-x_1$ in $f_1$ by $1-a_1x_1^{\delta_1},a_1>0,\delta_1>0,1-a_1x_1^{\delta_1}>0$ and by replacing the exponent $ax_2$ in $f_2$ by $ax_2^{\delta_2},\delta_2>0$. Then, one will end up with a H-function and several results will be available on this H-function. It may be observed that if in $f_1$, $\alpha$ is replaced by $\alpha-1$ and if $f_2$ is a general function then the first integral representation in Theorem 3.5 is also Erd\'elyi-Kober fractional integral of the first kind of order $\beta$ and parameter $\alpha$, see also Mathai and Haubold (2018).

 \vskip.3cm\noindent{\bf Problem 3.6.\hskip.3cm Type-1 beta versus type-1 beta}

 \vskip.3cm Let
 $$f_j(x_j)=\frac{\Gamma(\alpha_j+1+\beta_j)}{\Gamma(\alpha_j+1)\Gamma(\beta_j)}x_j^{\alpha_j}(1-x_j)^{\beta_j-1},0\le x_j\le 1,\alpha_j>-1,\beta_j>0,j=1,2.
 $$Let the joint function be $f_1(x_1)f_2(x_2)$ and let $u=\frac{x_2}{x_1}$. Then, proceeding as in the derivation of Problem 3.1, we have the following results:

 \vskip.3cm\noindent{\bf Theorem 3.6.}\hskip.3cm{\it For } $\alpha_j>-1,\beta_j>0,j=1,2$
 \begin{align*}
 G&=G^{1,1}_{2,2}\left[u\vert_{\alpha_2,-1-\alpha_1-\beta_1}^{-1-\alpha_1,\alpha_2+\beta_2}\right]\\
 &=\begin{cases}\frac{1}{\Gamma(\beta_1)\Gamma(\beta_2)}u^{\alpha_2}\int_0^{\frac{1}{u}}v^{\alpha_1+\alpha_2+1}(1-v)^{\beta_1-1}
 (1-uv)^{\beta_2-1}{\rm d}v,0\le u<1\\
 \frac{1}{\Gamma(\beta_1)\Gamma(\beta_2)}u^{-\alpha_1-\beta_1-1}\int_0^1v^{\alpha_1+\alpha_2+1}(u-v)^{\beta_1-1}(1-v)^{\beta_2-1}{\rm d}v,1\le u<\infty\end{cases}\end{align*}

 \begin{align*}
 G&=\begin{cases}\frac{\Gamma(2+\alpha_1+\alpha_2)}{\Gamma(2+\alpha_1+\alpha_2+\beta_1)\Gamma(\beta_2)}u^{\alpha_2}{_2F_1}(1-\beta_2,2+\alpha_1+\alpha_2;2+\alpha_1+\alpha_2+\beta_1;u),0\le u<1\\
 \frac{\Gamma(2+\alpha_1+\alpha_2)}{\Gamma(2+\alpha_1+\alpha_2+\beta_2)\Gamma(\beta_1)}u^{-2-\alpha_1}{_2F_1}(1-\beta_1,2+\alpha_1+\alpha_2;2+\alpha_1+\alpha_2+\beta_2;\frac{1}{u}),\\
 \mbox{for } 1\le u<\infty.\end{cases}\end{align*}

 \vskip.3cm A direct generalization will be available by replacing $1-x_j$ by $1-a_jx_j^{\delta_j},a_j>0,\delta_j>0,1-a_jx_j^{\delta_j}>0, j=1,2.$ Then, we will end up with a H-function and several results will be available on this H-function. One general case of the ratio $u=\frac{x_2}{x_1}$ will be discussed here for the sake of illustration. This will be listed here as Problem 3.7.

 \vskip.3cm\noindent{\bf Problem 3.7.\hskip.3cm Generalized type-2 beta versus generalized gamma}
 \vskip.3cm
 Let
 $$f_1(x_1)=\delta\frac{a^{\frac{\alpha+1}{\delta}}\Gamma(\frac{\alpha+1}{\delta}+\beta)}{\Gamma(\frac{\alpha+1}{\delta})\Gamma(\beta)}x_1^{\alpha}(1+ax_1^{\delta})^{-(\frac{\alpha+1}{\delta}+\beta)},f_2(x_2)=\frac{\rho b^{\frac{\gamma}{\rho}}}{\Gamma(\frac{\gamma+1}{\rho})}x_2^{\rho}{\rm e}^{-bx_2^{\rho}}
 $$for $\alpha>-1,\gamma>-1,\delta>0,\rho>0, a>0, b>0, \beta>0$ and $f_1$ and $f_2$ are zero elsewhere. Let $u=\frac{x_2}{x_1}$ with the associated function $g_2(u)$. The Mellin transforms of $f_1$ and $f_2$ are the following:
 $$M_{f_2}(s)=\frac{b^{\frac{1}{\rho}}}{\Gamma(\frac{\gamma+1}{\rho})}\Gamma(\frac{\gamma+s}{\rho})(b^{\frac{1}{\rho}})^{-s},\Re(\frac{\gamma+s}{\rho})>0$$and
 $$M_{f_1}(2-s)=\frac{1}{a^{\frac{1}{\delta}}\Gamma(\frac{\alpha+1}{\delta})\Gamma(\beta)}\Gamma(\frac{2+\alpha}{\delta}-\frac{s}{\delta})\Gamma(\beta-\frac{1}{\delta}+\frac{s}{\delta})(a^{\frac{1}{\delta}})^{s}
 $$for $\Re(2+\alpha-s)>0, \Re(\beta\delta -1+s)>0$. Then
 \begin{align*}
 M_{f_1}(2-s)M_{f_2}(s)&=\frac{b^{\frac{1}{\rho}}}{a^{\frac{1}{\delta}}\Gamma(\frac{\alpha+1}{\delta})\Gamma(\frac{\gamma+1}{\rho})\Gamma(\beta)}\\
 &\times \Gamma(\frac{\gamma}{\rho}+\frac{s}{\rho})\Gamma(\beta-\frac{1}{\delta}+\frac{s}{\delta})\Gamma(\frac{2+\alpha}{\delta}-\frac{s}{\delta})\left(\frac{b^{\frac{1}{\rho}}}{a^{\frac{1}{\delta}}}\right)^{-s}.\end{align*}Let
 $$c=\frac{\delta\rho a^{\frac{\alpha+1}{\delta}}b^{\frac{\gamma+1}{\rho}}\Gamma(\frac{\alpha+1}{\delta}+\beta)}{\Gamma(\frac{\alpha+1}{\delta})\Gamma(\beta)\Gamma(\frac{\gamma+1}{\rho})},c_1=\frac{b^{\frac{1}{\rho}}}{a^{\frac{1}{\delta}}
 \Gamma(\frac{\alpha+1}{\delta})\Gamma(\beta)\Gamma(\frac{\gamma+1}{\rho})}.
 $$Then, from the inverse Mellin transform we have
 \begin{align*}
 g_2(u)&=\frac{1}{2\pi i}\int_{c-i\infty}^{c+i\infty}M_{f_1}(2-s)M_{f_2}(s)u^{-s}{\rm d}s\\
 &=c_1\frac{1}{2\pi i}\int_{c-i\infty}^{c+i\infty}\Gamma(\frac{\gamma}{\rho}+\frac{s}{\rho})\Gamma(\beta-\frac{1}{\delta}+\frac{s}{\delta})\Gamma(\frac{2+\alpha}{\delta}-\frac{s}{\delta})
 \left[\frac{b^{\frac{1}{\rho}}u}{a^{\frac{1}{\delta}}}\right]^{-s}{\rm d}s\\
 &=c_1H_{1,2}^{2,1}\left[\frac{b^{\frac{1}{\rho}}u}{a^{\frac{1}{\delta}}}\big\vert_{(\frac{\gamma}{\rho},\frac{1}{\rho}),(\beta-\frac{1}{\delta},\frac{1}{\delta})}^{(1-\frac{2+\alpha}{\delta},\frac{1}{\delta})}\right],0\le u<\infty.\end{align*}From the integral representations of $g_2(u)$ we have the following forms and for $\delta=\rho$ and for $\gamma+\rho\nu\ne\beta\delta+\delta\lambda,v=0,1,...,\lambda=0,1,...$ the poles of the integrand are simple and then we have the following series representations for the H-function for $\delta=\rho$. We will state this as a theorem.

 \vskip.3cm\noindent{\bf Theorem 3.7.}\hskip.3cm{\it For the conditions, $c,c_1$ stated above, we have
 \begin{align*}
 g_2(u)&=cu^{-\alpha-2}\int_0^{\infty}v^{\alpha+\gamma+1}[1+a(\frac{v}{u})^{\delta}]^{-(\frac{\alpha+1}{\delta}+\beta)}{\rm e}^{-bv^{\rho}}{\rm d}v\\
 &=cu^{-\gamma}\int_0^{\infty}v^{\alpha+\gamma+1}[1+av^{\delta}]^{-(\frac{\alpha+1}{\delta}+\beta)}{\rm e}^{-b(uv)^{\rho}}{\rm d}v\\
 &=c_1\mbox{ H-function above}\end{align*}where for $\rho=\delta$,}
 \begin{align*}
 H&=H^{2,1}_{1,2}\left[(\frac{b}{a})^{\frac{1}{\delta}}u\big\vert_{(\frac{\gamma}{\delta},\frac{1}{\delta}),(\beta-\frac{1}{\delta},\frac{1}{\delta})}^{(1-\frac{2+\alpha}{\delta},\frac{1}{\delta})}\right]\\
 &=\delta\begin{cases}\Gamma(\beta-\frac{\gamma+1}{\delta})\Gamma(\frac{2+\alpha+\gamma}{\delta})[(\frac{b}{a})^{\frac{1}{\delta}}u]^{\gamma}{_1F_1}(\frac{2+\alpha+\gamma}{\delta};1+\frac{\gamma+1}{\delta}-\beta;\frac{bu^{\delta}}{a})\\
 +\Gamma(\frac{\gamma+1}{\delta}-\beta)\Gamma(\frac{\alpha+1}{\delta}+\beta)\left[(\frac{b}{a})^{\frac{1}{\delta}}u\right]^{\beta\delta-1}{_1F_1}(\frac{\alpha+1}{\delta}+\beta;1+\beta-\frac{\alpha+1}{\delta};\frac{b}{a}u^{\delta}),\\
 \mbox{ for }0\le u<\infty, \gamma-\delta\beta+1\ne\pm\nu,\nu=0,1,...\end{cases}\end{align*}

 \vskip.3cm\noindent{\bf 4.\hskip.3cm Generalization to Matrix-variate Cases}
 \vskip.3cm
Here, all the matrices appearing are $p\times p$ real symmetric positive definite when in the real domain or $p\times p$ Hermitian positive definite when in the complex domain. Let $X_1>O$ and $X_2>O$ be $p\times p$ real symmetric positive definite matrices with the associated functions $f_1(X_1)$ and $f_2(X_2)$ respectively and with the joint function $f_1(X_1)f_2(X_2)$, capital letters representing matrices and lower case letters denoting scalar variables where $f_j(X_j)$ is a real-valued scalar function of $X_j,j=1,2$. When $f_j(X_j)\ge 0$ for all $X_j$, in the domain of $X_j$, with $\int_{X_j}f_j(X_j){\rm d}X_j=1$ then $f_j(X_j)$ is a statistical density also. Since we are dealing with symmetric matrices we will be defining symmetric product corresponding to the scalar case $x_1x_2$ as $X_2^{\frac{1}{2}}X_1X_2^{\frac{1}{2}}$ where $X_2^{\frac{1}{2}}$ is the symmetric positive definite square root of $X_2$, and symmetric ratio corresponding to $\frac{x_2}{x_1}$ as $X_2^{\frac{1}{2}}X_1^{-1}X_2^{\frac{1}{2}}$ where $X_1^{-1}$ is the regular inverse of $X_1$. We can also define symmetric product corresponding to $x_1x_2$ as $X_1^{\frac{1}{2}}X_2X_1^{\frac{1}{2}}$ and symmetric ratio corresponding to $\frac{x_2}{x_1}$ as $X_1^{-\frac{1}{2}}X_2X_1^{-\frac{1}{2}}$ also.

\vskip.3cm\noindent{\bf 4.1.\hskip.3cm Symmetric product in the real matrix-variate case}
\vskip.3cm
Let $U=X_2^{\frac{1}{2}}X_1X_2^{\frac{1}{2}}$ with $V=X_2$ so that $X_1=V^{-\frac{1}{2}}UV^{-\frac{1}{2}}$. We can show that
$${\rm d}X_1\wedge{\rm d}X_2=|V|^{-\frac{p+1}{2}}{\rm d}U\wedge{\rm d}V\eqno(4.1)
$$see, for example Mathai (1997), where $|(\cdot)|$ denotes the determinant of $(\cdot)$. Let $X_j=(x_{jik})$ with the $(i,k)$-th element $x_{jik}$ for $i=1,...,p,k=1,...,p$. Since $X_j=X_j'$, a prime denoting the transpose, there are only $p(p+1)/2$ distinct elements in $X_j$ and only $p(p+1)/2$ differentials ${\rm d}x_{jik},i\le k$ or $i\ge k$. Then, ${\rm d}X_j=\wedge_{i\le k}{\rm d}x_{jik}$. Note that the wedge product in (4.1) remains the same if we take $V=X_1$ instead of $X_2$. For illustrative purposes, we will discuss one problem involving M-convolution of a product. For M-convolutions and M-transforms, see Mathai (1997).

\vskip.3cm\noindent{\bf Problem 4.1.\hskip.3cm Real $p\times p$ matrix-variate gamma versus  real $p\times p$ matrix-variate gamma}
\vskip.3cm
Let
$$f_j(X_j)=\frac{|B_j|^{\alpha_j}}{\Gamma_p(\alpha_j)}|X_j|^{\alpha_j-\frac{p+1}{2}}{\rm e}^{-{\rm tr}(B_jX_j)},X_j>O,B_j>O,j=1,2
$$where $X_j=X_j'>O$ is a real matrix of a $p(p+1)/2$ distinct (functionally independent) real scalar variables $x_{jik}$ as elements, $B_j=B_j'>O$ is a $p\times p$ constant positive definite matrix, ${\rm tr}(\cdot)$ means the trace of $(\cdot)$ and $\Gamma_p(\alpha_j)$ is a real matrix-variate gamma defined as
\begin{align*}
\Gamma_p(\alpha)&=\pi^{\frac{p(p-1)}{4}}\Gamma(\alpha)\Gamma(\alpha-\frac{1}{2})...\Gamma(\alpha-\frac{p-1}{2}),\Re(\alpha)>\frac{p-1}{2}\tag{4.2}\\
&=\int_{Y>O}|Y|^{\alpha-\frac{p+1}{2}}{\rm e}^{-{\rm tr}(Y)}{\rm d}Y,Y>O,\Re(\alpha)>\frac{p-1}{2}\tag{4.3}\end{align*}
and thus, the real matrix-variate gamma $\Gamma_p(\alpha)$ is associated with the real matrix-variate gamma integral in (4.3). Hence, $\Gamma_p(\alpha)$ is called as the real matrix-variate gamma. Let $g_2(U)$  be the real-valued scalar function associated with the matrix $U=X_2^{\frac{1}{2}}X_1X_2^{\frac{1}{2}}$. Then, from the joint function $f_1(X_1)f_2(X_2)$ we have for $c=\prod_{j=1}^2\frac{|B_j|^{\alpha_j}}{\Gamma_p(\alpha_j)}$,
\begin{align*}
g_2(U)&=c\int_{V>O}|V|^{-\frac{p+1}{2}}|V^{-\frac{1}{2}}UV^{-\frac{1}{2}}|^{\alpha_1-\frac{p+1}{2}}|V|^{\alpha_2-\frac{p+1}{2}}{\rm e}^{-{\rm tr}(B_1V^{-\frac{1}{2}}UV^{-\frac{1}{2}})-{\rm tr}(B_2V)}{\rm d}V\\
&=c|U|^{\alpha_1-\frac{p+1}{2}}\int_{V>O}|V|^{\alpha_2-\alpha_1-\frac{p+1}{2}}{\rm e}^{-{\rm tr}(B_1V^{-\frac{1}{2}}UV^{-\frac{1}{2}})-{\rm tr}(B_2V)}{\rm d}V\tag{4.4}\\
&=c|U|^{\alpha_2-\frac{p+1}{2}}\int_{V>O}|V|^{\alpha_1-\alpha_2-\frac{p+1}{2}}{\rm e}^{-{\rm tr}(B_1V)-{\rm tr}(B_2V^{-\frac{1}{2}}UV^{-\frac{1}{2}})}{\rm d}V.\tag{4.5}\end{align*}
Note that (4.4) and (4.5) can be taken as real matrix-variate generalization of Bessel integrals, Kr\"atzel integral, Bayesian structures etc. The integrand, normalized can also be taken as a real matrix-variate inverse Gaussian density when $f_j(X_j),j=1,2$ are densities. In this case, $g_2(U)$ is the unique density of the symmetric product $U$ also. All the pairs of functions that we considered in Section 2 can be generalized to the corresponding matrix-variate cases except Problem 2.7. In Problem 2.7, generalization to the matrix-variate case is possible for some specific values of $\delta_1$ and $\delta_2$.

\vskip.3cm\noindent{\bf Concluding Remarks}
\vskip.3cm
In Sections 2 and 3, elementary functions belonging to the pathway family are considered so that the results in Sections 2 and 3 will be readily applicable to various practical problems. Immediate generalizations of all pairs in Sections 2 and 3 are also useful in situations in real scalar variable cases. All problems in Sections 2 and 3, except Problems 2.7 and 3.7, can be generalized to the corresponding matrix-variate cases in the real and complex domains. In all those cases, when $f_j(X_j),j=1,2$ are statistical densities then $g_1(U)$ and $g_2(U)$ will be unique densities corresponding to symmetric product and symmetric ratios respectively. There are  not many results in the literature on Mellin convolutions of products and ratios involving more than two real scalar variables. Let $u=x_1x_2...x_k$ with the joint function $f_1(x_1)...f_k(x_k)$. Then, the integral representations and series representations of $g_1(u)$ can be obtained in many different ways which will produce various results. When it comes to the ratio, then ratios can be defined in many different ways when more than two real scalar variables are involved. Then, the corresponding $g_2(u)$ can take different forms and shapes.
\vskip.2cm
``Convolution'' in mathematical literature often means Laplace convolution. Mellin convolutions of products and ratios used to appear very rarely in mathematical literature. In the area of special functions, when the functions are defined through Mellin-Barnes representations, Mellin transforms and Mellin convolutions appear automatically. The class of hypergeometric functions appear in very many practical problems as models and as series solutions of differential equations. Hypergeometric function ${_pF_q}(\cdot)$ has a Mellin-Barnes representation. Mittag-Leffler function is considered as the queen function in the area of fractional calculus. These functions also have Mellin-Barnes representations. Wright's function is another prominent function in the literature having Mellin-Barnes representation. Generalized scalar variable special functions such as G and H-functions have Mellin-Barnes representations. In all these cases, Mellin convolutions of products and ratios appear. Products and ratios of random variables appear in many real-life situations. In industrial production processes the money value of the output has the structure of the product of two real scalar positive random variables, namely, quantity of output and price per unit. In a series of papers, starting from 2009, Mathai has shown that all types fractional integrals introduced by various authors from time to time are nothing but Mellin convolutions of products and ratios of real scalar variables. In the overview paper Mathai (2024) it is shown that the whole area of fractional calculus, including fractional integrals, fractional derivatives etc can be studied by using Mellin convolutions of products and ratios. It is shown that the approach through Mellin convolutions enable us to go to the Mellin transforms thereby to the corresponding inverse Mellin transforms which provide explicit evaluations of fractional integrals. It is also shown that when the arbitrary functions in fractional integrals belong to some general classes, then the inverse Mellin transform  will end up in G and H-functions. Then, one could make use of the G-function differential equation to construct differential equations for fractional integrals. It is also shown in Mathai (2024) that the approach through statistical distribution theory or Mellin convolutions enable us to extend fractional integrals to functions of matrix arguments and to the complex domains.
\vskip.2cm
Explicit computable series representations for various categories of G-functions, including the general G-function, are given in Mathai (1993) book. The same techniques work in obtaining explicit computable series forms for the H-function also. With the help of these computable representations, computer packages are also produced. Programs are available in MAPLE and MATHEMATICA for analytical and numerical computations of G and H-functions Hence, numerical computations and graphs are not attempted in the present paper.
\vskip.2cm
In the case of Mellin convolutions of products involving two functions, it is shown that when these two functions belong to generalized gamma densities then one can obtain Kr\"atzel integrals and Kr\"atzel transforms, reaction-rate probability integrals in nuclear reaction-rate theory and inverse Gaussian density in stochastic processes. The large contributions of Mathai and Haubold in astrophysics area started with the evaluation of a reaction-rate probability integral by using the density of product of real scalar positive random variables, which is nothing but the Mellin convolution of a product involving generalized gamma densities. A summary of the work in this area of astrophysics is available in the monograph Mathai and Haubold (1988) and an updated version of this book is scheduled to appear in early 2025 from Springer.
\vskip.2cm
In the matrix-variate case, Mathai (1997) defined M-convolutions and M-transforms, corresponding to Mellin convolutions and Mellin transforms in the real scalar case. Jacobians of matrix transformations and functions of matrix arguments in the real and complex domains are given there, along with the necessary conditions for the existence of various items. Detailed existence conditions for G-functions are given in Mathai (1993). These M-convolutions are shown to be densities of symmetric product and symmetric ratios of matrices when the basic functions involved are matrix-variate statistical densities in the real or complex domain. Distributions of symmetric products and symmetric ratios of real positive definite or Hermitian positive definite matrix-variate random variables appear in a large number of practical situations such as in neutrino problems [Fortin, Giasson and Marleau (2016)], in matrix-variate fractional calculus [Mathai (2024)], in the analysis of multiple and multi-look data in radar, sonar etc [Deng (2016)]. Matrix-variate statistical distributions including singular distributions in the real and complex domains and their applications in various areas are given in Mathai, Provost and Haubold (2022). In this book, detailed existence conditions of the various results in the real and complex domains are also given. 

\vskip.3cm\noindent{\begin{center}{\bf References}\end{center}

\vskip.2cm\noindent X. Deng (2016: {\it Texture Analysis and Physical Interpretation of Polarimetric SAR Data}, Ph.D thesis, Universitat Politecnica de Catelunya, Barcelona, Spain, 2016.
\vskip.2cm\noindent J.-F. Fortin, N. Giasson, and L. Marleau (2016): Probability density function for neutrino masses and mixings, {\it Physical Review} {\bf D 94}, 115004.
\vskip.2cm\noindent J.-F. Fortin, N. Giasson, L. Marleau, and J. Pelletier-Dumont (2020): Mellin transform approach to rephasing invariants, {\it Physical Review} {\bf D 102}, 036001.
\vskip.2cm\noindent Y. Luchko (2008): Integral transforms of the Mellin convolution type and their generating operators, {\it Integral Transforms and Special Functions}, {\bf 19}, 809-851.
\vskip.2cm\noindent A.M. Mathai (1993): {\it  A Handbook of Generalized Special Functions for Statistical and Physical Sciences}, Oxford University Press, Oxford.
\vskip.2cm\noindent A.M. Mathai (1997): {\it Jacobians of Matrix Transformations and Functions of Matrix Arguments}, World Scientific Publishing, New York.
\vskip.2cm\noindent A.M. Mathai (2005): A pathway to matrix-variate gamma and normal densities, {\it Linear Algebra and its Applications}, {\bf 396}, 317-328.
\vskip.2cm\noindent A.M. Mathai (2024): Statistical distribution theory and fractional calculus, {\it Stats2024}, {\bf 7}, 1259-1295. http:// doi.org/ 10.3390/ stats7040074.
\vskip.2cm\noindent A.M. Mathai and Hans J. Haubold (1988): {\it Modern Problems in Nuclear and Neutrino Astrophysics}, Akademie-Verlag, Berlin, 1988 (Updated version from Spinger Nature in 2025).
\vskip.2cm\noindent A.M. Mathai and H.J. Haubold (2018): {\it Erd\'elyi-Kober Fractional Calculus from a Statistical Perspective, Inspired by Solar Neutrino Physics}, Springer Briefs in Mathematical Physics 31, Springer Nature, Singapore.
\vskip.2cm\noindent A.M. Mathai and H.J. Haubold (2020): Mathematical aspects of Kr\"atzel integral and Kr\"atzel transform, {\it mathematics}, {\bf 8}, 526; doi:10.3390/math8040526.
\vskip.2cm\noindent A.M. Mathai and S.B. Provost (2006): Some complex matrix-variate statistical distributions on rectangular matrices, {\it Linear Algebra and its Applications}, {\bf 40}, 198-216.
\vskip.2cm\noindent Arak M. Mathai, Serge B. Provost and Hans J. Haubold (2022): {\it Multivariate Statistical Analysis in the Real and Complex Domains}, Springer Nature, Switzerland, 2022.
\vskip.2cm\noindent A.M. Mathai, R.K. Saxena and H.J. Haubold (2010): {\it The H-function: Theory and Applications}, Springer, New York.
\vskip.2cm\noindent F. Mainardi and G. Pagnini (2008): Mellin-Barnes integrals for stable distributions and their convolutions, {\it Fractional Calculus and Applied Analysis}, {\bf 11}, 443-456.

 \end{document}